\documentclass[natbib]{svmult}

\usepackage{tikz} 
\usepackage[tight]{subfigure}
\usepackage{mathptmx}      
\usepackage{textcomp}
\usepackage{rotating}
\usepackage{array}
\usepackage{booktabs}
\usepackage{enumitem}
\usepackage{mathbbol}

\usetikzlibrary{arrows.meta,shapes.geometric,trees}
\usepackage{tikz-qtree}
\renewcommand{\cite}{\citep}

\newtheorem{scheme}{Argumentation Scheme}

\title{Argumentation in Mathematical Practice}
\author{Andrew Aberdein
\and Zoe Ashton}
\institute{A. Aberdein \at
              School of Arts and Communication, Florida Institute of Technology,
150 West University Blvd,
Melbourne, FL 32901, USA \\
              \email{aberdein@fit.edu}           
           \and
Z. Ashton \at Department of Philosophy, Ohio State University, 230 N Oval Mall, Columbus, OH 43210, USA}

\begin{document}
	
\maketitle

\abstract{
Formal logic has often been seen as uniquely placed to analyse mathematical argumentation. While formal logic is certainly necessary for a complete understanding of mathematical practice, it is not sufficient. Important aspects of mathematical reasoning closely resemble patterns of reasoning in non-mathematical domains. Hence the tools developed to understand informal reasoning, 
collectively known as argumentation theory, are also applicable to much mathematical argumentation. 
This chapter investigates some of the details of that application. Consideration is given to the many contrasting meanings of the word `argument'; to some of the specific argumentation-theoretic tools that have been applied to mathematics, notably Toulmin layouts and argumentation schemes; to some of the different ways that argumentation is implicated in mathematical practices; and to the social aspects of mathematical argumentation.
}

\section{Introduction}\label{section:Intro}

Since logic developed the tools to adequately represent formal derivations, many philosophers of mathematics have been tempted to conclude that formal derivation suffices to account for all interesting features of mathematical practice.
However, there have always been other philosophers who perceived the shortcomings of such a reduction. Here, for example, is Henri Poincar\'e:
\begin{quote}
If you are present at a game of chess, it will not suffice, for the understanding of the game, to know the rules for moving the pieces. That will only enable you to recognize that each move has been made conformably to these rules, and this knowledge will truly have very little value. Yet this is what the reader of a book on mathematics would do if he were a logician only. To understand the game is wholly another matter; it is to know why the player moves this piece rather than that other which he could have moved without breaking the rules of the game. It is to perceive the inward reason which makes of this series of successive moves a sort of organized whole. This faculty is still more necessary for the player himself, that is, for the inventor \cite[218]{Poincare13}.
\end{quote}
One response to this limitation of formal logic is to recognise an analogy with a similar limitation in another domain: formal logic is also an imperfect tool for understanding everyday reasoning.
Solutions have been proposed for that problem: systems of informal logic, argumentation theory, or dialectic have been devised since antiquity to address ordinary reasoning \cite{Eemeren14}. 
Hence some philosophers of mathematical practice have reasoned that these theories might also lend themselves to the understanding of mathematical reasoning \cite{Aberdein09b}.

This chapter surveys the uses to which argumentation theory has been put in order to understand mathematical practice.
Section~\ref{section:Argument} addresses the many ambiguities implicit in the word `argument'---and how they are specifically related to mathematics. 
Section~\ref{section:Proof} discusses two prominent proposals for the application of particular tools from argumentation theory to mathematics: Toulmin layouts and argumentation schemes.
Section~\ref{section:Reasoning} considers the argumentative aspects of mathematical practices beyond proof and section~\ref{section:Communities} 
focuses on the contribution of communities of mathematical practice to such argumentation.

\section{What is an argument?}\label{section:Argument}

Several overlapping distinctions may be drawn in our understanding of arguments. 

\begin{itemize}
\item {\bf Argument-that/argument-about}
A sequence of statements, whereby premises offer support for a conclusion, is an argument. But an exchange of conflicting views held by different people is also an argument. We may refer to the former as an \emph{argument-that} and to the latter as an \emph{argument-about}. 
Arguments-that are also known as arguments$_1$ and arguments-about as arguments$_2$ \citep{OKeefe77}.
As Michael Gilbert helpfully glosses the distinction, ``one person makes an argument$_1$ and [at least] two people have an argument$_2$'' \cite[21]{Gilbert14}.
Although the most familiar mathematical arguments tend to be arguments-that, arguments-about 
arise in mathematics too, such as priority disputes, contested axioms or principles, or debates over the legitimacy of a technique or the admissibility of a proof. A salient recent example is the contested status of Shinichi Mochizuki's claimed proof of the $abc$ conjecture \cite{Aberdein23}.

\item {\bf Process/product} Arguments-that are sometimes represented as products of the argument-about process, but this is arguably a mischaracterization: 
``If, as part of organizing the domain of argumentation theory, we merely want to distinguish acts of arguing from arguments-as-objects, we should not use the misleading process/product labels to do so. At the very least such labels imply a relationship that does not exist and so distort our perceptions of the domain of study
'' \cite[87]{Goddu11}.
Nonetheless, we may usefully distinguish between an argument understood as an act of arguing and the informational trace that act leaves behind (a transcript, a recording,\dots), also often called an argument \cite[948]{Sundholm12a}.
This distinction 
straightforwardly applies to mathematics.
(For further reflection on proofs as acts of proving or proof-events, see \citealp{Goguen01,Stefaneas14,Stefaneas15}.)

\item {\bf Monologue/dialogue/polylogue} The number of participants in an argument may vary considerably. Most attention has traditionally been paid to monologues and dialogues: arguments-that are characteristically presented as monologues; arguments-about as dialogues. If dialogues are understood loosely, as also covering argumentation with more than two participants, that distinction would exhaust the options. However, some authors have made a case for differentiating the two-participant dialogue from the many-participant polylogue \citep{Lewinski14}.
The Polymath Project, a series of crowd-sourced proofs of open conjectures, is a rich source for research on 
mathematical 
polylogues \cite{Allo21}. 

\item {\bf Small scale/large scale} Arguments can vary significantly in scale. The scale of an argument may be measured in several conceptually distinct 
ways. The duration of an argument 
may range from arguments 
that take seconds to arguments which last for hundreds of years. 
The size of an argument 
may range from a few inferences expressed in short sentences about a simple issue to inferentially complex structures involving 
a great many very long sentences.
Likewise, mathematical proofs can vary in length from a few lines to tens or even hundreds of pages. In exceptional cases, proofs can be so long as to defy the capacity of any single mathematician to survey the whole \cite{Coleman09}.
This can arise in proofs achieved by traditional means, such as that of the classification of finite simple groups, the components of which comprised thousands of pages in several hundred articles by dozens of mathematicians \cite{Steingart12}.
It is even more acute in the case of computer assisted proofs, such as that of the four colour theorem or the Kepler conjecture, or indeed the subsequent computer verification of these proofs \cite{Gonthier08a,Hales17}.

\item {\bf Static/dynamic} Static arguments have achieved a final and definite form; dynamic arguments are fluid and ongoing. In general, the evolution of knowledge may be understood as the product of dynamic argumentation. Dynamic arguments 
are common in everyday life---but they are also 
central to the development of scientific thought. 
A well-known and influential analysis of a dynamic argument in mathematics is Imre Lakatos's \emph{Proofs and Refutations}. Lakatos shows how the protracted search for a proof of the Descartes--Euler conjecture, which relates the quantities of vertices, edges, and faces of convex polyhedra, involved significant redefinition of most of the concepts used in that conjecture \cite{Lakatos76}.
(For further discussion, see Section~\ref{subsection:Refutations} below or \citealp{Reyes21}.)

\item {\bf Centralized/distributed} Arguments may be either centralized or distributed with respect to several factors including time, people, space, and media. A highly centralized argument may be restricted to what one person communicates in one place at one time in one mode of expression. But arguments may involve a varying cast of arguers and be drawn out over long periods, in multiple locations and media. 
Large scale arguments are characteristically distributed: as a mathematical example, the classification of finite simple groups involved many years of work by a large, geographically widespread collective of mathematicians \cite{Steingart12}.

\item {\bf Sequential/parallel} 
Logic, whether formal or informal, tends to reconstruct arguments as a linear sequence of premises from which intermediate statements are derived, culminating in a final conclusion. (Some logical systems prefer to invert this sequence.) However, argumentation in natural contexts often occurs in a more parallel fashion: several strands of argument may be developed simultaneously, the conclusion may be derived from initial premises and then reinforced by subsequent subarguments, and so forth.
Parallel arguments are also more likely to arise in projects with many participants, whether by accident or design:
``To permit a large collaboration, \dots\ 
long proofs have been broken up into series of shorter lemmas'' \cite[11]{Hales17}.
The educationists Christine Knipping and David Reid have proposed more fine-grained subdivisions of parallel arguments in mathematics (\citealp{Knipping19}; for further discussion, see Section~\ref{subsection:Toulmin} below.)
\end{itemize}

The seven overlapping distinctions addressed so far have received unequal interest in the philosophy of mathematics. The formal logical approach mentioned in the introduction best coheres with a small scale, static, centralized, sequential product like a published proof. Much recent work in the philosophy of mathematical practice focuses on mathematical arguments along the other dimensions. 
Some further distinctions arise from consideration of the 
goals of the arguers:

\begin{itemize}
\item {\bf Persuasive/directive/polemic/\dots}
That arguments may be distinguished by their objective is an ancient idea: Aristotle distinguished forensic arguments (concerned with past acts), display arguments (concerned with present circumstances), and deliberative arguments (concerned with future acts) \citep[1358b]{Aristotle91}. 
Erik C. W. Krabbe and Jan Albert van Laar propose an updated distinction between three different functions of reasoning: persuasive (to convince the other party), directive (to get the other party to act), and polemic (to intimidate the other party) \cite[29 f.]{Krabbe07}.
They contrast these ``inherently argumentative'' functions with three further functions of reasoning: explanatory (to enhance understanding), explorative (to investigate connections between statements), and probative (to establish new knowledge).

\item {\bf Adversarial/non-adversarial}
On a strict interpretation, all arguments begin in conflict: a difference in belief, or concerning how to act, or of some other kind. 
On a broader interpretation, arguments need not be strictly adversarial, hence they may proceed from other situations, including shared uncertainty or one party knowing what another does not.
\end{itemize}

These nine different dimensions of comparison interact in important ways. Firstly, 
they are not pairwise independent. For example, as the scale of an argument increases, we may tentatively expect both the number of participants and the number of different objectives 
to increase, but the likelihood of the argument being either static, centralized, or sequential to tend to zero.
Secondly, the last two dimensions combine to produce 
what Douglas Walton calls \emph{dialogue types}
(see Fig.~\ref{dialoguetree}, adapted from \citealp[81]{Walton95}). 
Further complicating this picture, Walton and Krabbe observe that dialogues can shift from one type to another (for example, an inquiry might turn into a persuasion dialogue if one inquirer becomes an advocate for a particular result) or be embedded in a dialogue of a different type (a deliberation over which course of action to pursue might contain an inquiry into the merits of one action, say).
Walton maintains that arguments may arise in any dialogue type \cite{Walton98}; Krabbe is more conservative and regards argumentation as restricted to adversarial contexts (
the lefthand fork of Fig.~\ref{dialoguetree})
\cite[33]{Krabbe07}.

\begin{figure}[tbp]
\begin{center}
\scriptsize
\begin{tikzpicture} 
\tikzset{grow'=down,level distance=.7in}
\tikzset{every tree node/.style={align=center,anchor=north}}
\Tree [.{Is there a conflict?} 	
	[.\fbox{NO}\\{Is there a common problem to be solved?} 	[.\fbox{NO}\\\textsc{Information Seeking} ] 
					[.\fbox{YES}\\{Is this a theoretical problem?} 	[.\fbox{NO}\\\textsc{Deliberation} ] 
													[.\fbox{YES}\\\textsc{Inquiry} ] ] ]  
	[.\fbox{YES}\\{Is resolution the goal?} 	[.\fbox{NO}\\{Is settlement the goal?} 	[.\fbox{NO}\\\textsc{Eristics} ] 
														[.\fbox{YES}\\\textsc{Negotiation} ] ] 
								[.\fbox{YES}\\\textsc{Persuasion} ] ] ]
\end{tikzpicture}
\end{center}
\caption{Determining the type of dialogue \cite[after][81]{Walton95}}\label{dialoguetree}
\end{figure}

How do these distinctions apply to argumentation in mathematics?
The most widely discussed case is that of mathematical proof, which many authors have maintained is intrinsically argumentative. But even here, we may distinguish multiple distinct activities which give rise to arguments of different kinds. For example, Krabbe observes the following distinct stages:
\begin{quote}
\begin{enumerate}[noitemsep]
\item thinking up a proof to convince oneself of the truth of some theorem; 
\item thinking up a proof in dialogue with other people (inquiry dialogue; probative functions of reasoning); 
\item presenting a proof to one's fellow discussants in an inquiry dialogue (persuasion dialogue embedded in inquiry dialogue; persuasive and probative functions of reasoning); 
\item presenting a proof to other mathematicians, e.g. by publishing it in a journal (persuasion dialogue; persuasive and probative functions of reasoning) 
\item presenting a proof when teaching (information-seeking and persuasion dialogue; explanatory, persuasive, and probative functions of reasoning) 
\cite[457]{Krabbe08}.
\end{enumerate}
\end{quote}
This sequence is familiar from many mathematicians' descriptions of the proving process, although in actual examples some steps may be repeated as proof attempts come unstuck.

\begin{table}[t]
\caption{Some mathematical dialogue types}\label{mathtypes}
\begin{center}\small
\begin{tabular}{*{5}{m{.17\textwidth}@{~}}}
\toprule
\textbf{
\mbox{Dialogue} Type}&\textbf{Initial \mbox{Situation}}&\textbf{Main Goal}&\textbf{\mbox{Goal~of} \mbox{Proponent
}}&\textbf{\mbox{Goal~of} \mbox{Respondent
}}\\
\midrule
\raggedright 
Inquiry&Open-mindedness&\raggedright {Prove or} disprove \mbox{conjecture}& \raggedright Contribute to main goal& Obtain \mbox{knowledge}\\
\addlinespace
\raggedright 
Persuasion& \raggedright Difference of opinion&\raggedright Resolve difference {of opinion with} rigour& \raggedright Persuade respondent
&Persuade \mbox{proponent
}\\
\addlinespace
\raggedright 
Pedagogical Information-Seeking&\raggedright Respondent 
lacks information&\raggedright Transfer of knowledge&\raggedright Disseminate knowledge of results and methods&Obtain \mbox{knowledge}\\
\addlinespace
\raggedright 
Oracular Information-Seeking
& \raggedright 
Proponent lacks information& \raggedright Transfer of knowledge&Obtain \mbox{information}&{Inscrutable}\\
\addlinespace
\raggedright 
Deliberation
&Open-mindedness& \raggedright Reach a provisional conclusion& \raggedright Contribute to main goal& Obtain \mbox{warranted} belief\\
\addlinespace
\raggedright 
Negotiation
& \raggedright Difference of opinion& \raggedright Exchange resources for a provisional conclusion&\raggedright Contribute to main goal&Maximize \mbox{value of} \mbox{exchange}\\
\addlinespace
\raggedright 
Eristic&\raggedright Personal conflict&
\raggedright Reveal deeper conflict&Win in the eyes of onlookers&Win in the eyes of onlookers\\
\bottomrule
\end{tabular}
\par\medskip
\end{center}
\end{table}

In deference to his intrinsically adversarial conception of argument, Krabbe only considers three dialogue types as hospitable to proofs.
In other work, one of us has suggested that proofs (or other mathematical arguments) may be found in other dialogue types too (see Table~\ref{mathtypes}, from \citealp[165]{Aberdein21a}; see also \citealp[148]{Aberdein07}).
``Oracular'' information-seeking owes its inspiration to an influential aside of Alan Turing concerning a machine ``supplied with some unspecified means of solving number-theoretic problems; a kind of oracle as it were'' \cite[172]{Turing39}.
An oracle is a ``black box''---it supplies answers but not explanations. For some sceptics of computer-assisted proofs, this is a compelling analogy for the role that the computer plays in such proofs \cite[e.g.][]{Tymoczko79}.
Deliberation differs from inquiry in seeking only a provisional conclusion. Mathematicians aspire to more permanent stability for their results. Nonetheless, there are circumstances where they are obliged to settle for less than they would wish, despite the rigour of their arguments. These include 
the ``architectural conjectures'' upon which many mathematical research programmes depend \cite[198]{Mazur97}.
Negotiation characteristically adds resource sensitivity to the provisional outcome typical of deliberation. 
While idealized accounts of mathematical practice disregard such factors, they are unavoidable in some contexts, especially in applied mathematics. (And it has been controversially suggested that ``semi-rigorous'' proofs might come with price tickets, costing the computational resources necessary for certainty \cite{Zeilberger93}.)
Even eristic dialogues can be the context for mathematical reasoning, as demonstrated in the mathematically inventive quarrels of early modern mathematicians such as Girolamo Cardano and Niccol\`o Tartaglia \cite{Toscano20}.

A final distinction amongst the different senses of argument cuts across most of those discussed above and arises not from argumentation theory but directly from 
mathematical practice. Here it is presented by Joel David Hamkins:
\begin{itemize}
\item {\bf Hard/soft}
``A hard argument is one that is technically difficult; perhaps it involves a laborious construction or a difficult calculation; perhaps it involves bringing disparate fine details together in just the right combination in order to succeed; or perhaps it involves proving various specific facts about a comparatively abstract construction, perhaps relating disparate levels of abstraction.
A soft argument, in contrast, is one that appeals only to very general abstract features of the situation, and one needs hardly to construct or compute anything at all''
\cite[166]{Hamkins20}.
\end{itemize}
Hamkins's use of `hard' and `soft' echoes G.~H.\ Hardy's division of analysis into ``the `hard, sharp, narrow' kind as opposed to the `soft, large, vague' kind'' \cite[64]{Hardy29}.
It also owes something to Alexander Grothendieck's celebrated analogy between two strategies for solving mathematical problems and two ways of opening a nut: cracking it with a hammer or softening it in water until it opens with light pressure \cite[301]{McLarty03}. Grothendieck favoured the latter strategy, of immersing problems in a much wider theory from which a solution could (eventually) be readily inferred.
This in turn is suggestive of Freeman Dyson's division of mathematicians into birds and frogs: ``Birds fly high in the air and survey broad vistas of mathematics out to the far horizon. They delight in concepts that unify our thinking and bring together diverse problems from different parts of the landscape. Frogs live in the mud below and see only the flowers that grow nearby. They delight in the details of particular objects, and they solve problems one at a time'' \cite[212]{Dyson09}. 
Hard arguments play to the tightly focused strengths of the frogs; soft arguments require the birds' sweeping perspective.

\section{Proof as Argumentation}\label{section:Proof}

The view mentioned in Section~\ref{section:Intro}, that formal derivation suffices to account for all interesting features of mathematical practice, squarely focuses on the role of mathematical proof. Pure mathematicians trade in proofs. But, even in the domain of mathematical proof, formal logic does not capture all there is of interest. Toulmin layouts and argumentation schemes, as we will see below, are useful methods of examining mathematical proofs. 
But, notably, such techniques tease apart the descriptive and normative aspects of (apparent) proof. Whereas a formal derivation is either correct or not, strictly speaking, a derivation at all, the methods of argumentation theory provide the means to describe proofs independently of whether they succeed as proofs.

\subsection{Toulmin layouts}\label{subsection:Toulmin}
\begin{figure}[t]
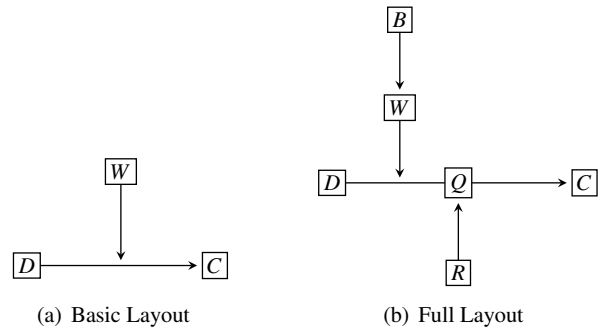

\label{Toulmin}
\begin{center}
\subfigure[Basic Layout]
{\label{ToulminA}\begin{pgfpicture}{0in}{0in}{1.5in}{1.5in}
\pgfsetlinewidth{.5pt} 
\pgfnodebox{NodeC1}[stroke]
   {\pgfxy(3.25,0.4)}
   {$C$}{2pt}{2pt} 
\pgfnodebox{NodeD1}[stroke]
   {\pgfrelative{\pgfxy(-2.5,0)}{\pgfnodecenter{NodeC1}}}
   {$D$}{2pt}{2pt}  
\pgfnodebox{Invisible1}[virtual]
   {\pgfrelative{\pgfxy(-1.25,0)}{\pgfnodecenter{NodeC1}}}{}{0pt}{0pt} 
\pgfnodebox{NodeW1}[stroke]
   {\pgfrelative{\pgfxy(0,1.25)}{\pgfnodecenter{Invisible1}}}
   {$W$}{2pt}{2pt}

\pgfsetarrowsend{stealth} 
\pgfnodesetsepstart{0pt}\pgfnodesetsepend{2pt} 
\pgfnodeconnline{NodeD1}{NodeC1} 
\pgfnodeconnline{NodeW1}{Invisible1} 
\end{pgfpicture}} 
\subfigure[Full Layout]
{\label{ToulminB}\begin{pgfpicture}{0in}{0in}{2in}{1.5in}
\pgfsetlinewidth{.5pt} 
\pgfnodebox{NodeQ1}[stroke]
   {\pgfxy(2.6,1.5)}
   {$Q$}{2pt}{2pt} 
\pgfnodebox{NodeD1}[stroke]
   {\pgfrelative{\pgfxy(-1.5,0)}{\pgfnodeborder{NodeQ1}{180}{0pt}}}
   {$D$}{2pt}{2pt}  
\pgfnodebox{NodeC1}[stroke]
   {\pgfrelative{\pgfxy(1.5,0)}{\pgfnodeborder{NodeQ1}{0}{0pt}}}
   {$C$}{2pt}{2pt}  
\pgfnodebox{NodeR1}[stroke]
   {\pgfrelative{\pgfxy(0,-1)}{\pgfnodeborder{NodeQ1}{270}{0pt}}}
   {$R$}{2pt}{2pt}  
\pgfnodebox{Invisible1}[virtual]
   {\pgfrelative{\pgfxy(-.6,0)}{\pgfnodeborder{NodeQ1}{180}{0pt}}}{}{0pt}{0pt} 
\pgfnodebox{NodeW1}[stroke]
   {\pgfrelative{\pgfxy(0,1)}{\pgfnodecenter{Invisible1}}}
   {$W$}{2pt}{2pt}
\pgfnodebox{NodeB1}[stroke]
   {\pgfrelative{\pgfxy(0,1)}{\pgfnodeborder{NodeW1}{90}{0pt}}}
   {$B$}{2pt}{2pt}

\pgfnodeconnline{NodeD1}{NodeQ1} 
\pgfsetarrowsend{stealth} 
\pgfnodesetsepstart{0pt}\pgfnodesetsepend{2pt} 
\pgfnodeconnline{NodeQ1}{NodeC1} 
\pgfnodeconnline{NodeR1}{NodeQ1} 
\pgfnodeconnline{NodeW1}{Invisible1} 
\pgfnodeconnline{NodeB1}{NodeW1} 
\end{pgfpicture}}
\caption{Toulmin layouts}
\end{center}
\end{figure}

One of the most influential treatments of informal argumentation is that of Stephen Toulmin \citeyearpar{Toulmin58}. 
His `layout' can represent deductive inference, but encompasses many other species of argument besides.  In its simplest form, shown in Fig.~\ref{ToulminA}, the layout represents the derivation of a Claim ($C$), from Data ($D$), in accordance with a Warrant ($W$).  This DWC pattern 
resembles a deductive inference rule, such as \textit{modus ponens}, but it can be used to represent looser inferential steps.  
The differences between the types of inference which the layout may represent are made explicit by the additional elements of the full layout shown in Fig.~\ref{ToulminB}.  
The warrant is justified by its dependence on Backing ($B$), possible exceptions or Rebuttals ($R$) are allowed for, and the resultant force of the argument is stated in the Qualifier ($Q$).  Hence the full layout may be understood as `{Given that} $D$, {we can}  $Q$ {claim that} $C$, {since} $W$ ({on account of} $B$), {unless} $R$'. In a frequently cited example (derived from \citealp[104]{Toulmin58}), `{Given that \textsc{Harry was born in Bermuda}
, we can \textsc{presumably} 
 claim that \textsc{he is British}
, since \textsc{anyone born in Bermuda will generally be British} 
  (on account of \textsc{various statutes~\dots}
), unless \textsc{his parents were aliens, say}
}.'

Toulmin wrote \emph{The Uses of Argument} in England in the 1950s as a critique of what he perceived as a formalizing trend in contemporary philosophy; the Toulmin layout was subsequently adopted by communication theorists in America in the 1960s and 70s; from there it seems to have passed to mathematics educationists in Germany in the 1990s. 
(Although Toulmin himself briefly considers a mathematical example, Theaetetus’s proof that there are exactly five platonic solids \cite[89]{Toulmin79}. For discussion, see \cite[290 ff.]{Aberdein05}.)
In particular, G\"otz Krummheuer is usually credited with the first sustained application of the Toulmin layout to mathematical argumentation (\citealp{Krummheuer95}; for a recent survey, see \citealp{Krummheuer15}). 

Toulmin draws a distinction between analytic and substantial arguments depending on whether the claim is already at least implicit in the backing.
This may suggest that only analytic arguments occur in mathematics. But that overlooks the emphasis Toulmin places on how out of the ordinary analytic arguments are: ``it begins to be a little doubtful whether any genuine, practical argument could ever be properly analytic'' \cite[127]{Toulmin58}. Granted, immediately after this passage, Toulmin places arguments in (pure) mathematics amongst the analytic arguments. However, the philosophers and educationists who have applied Toulmin's work to mathematics endorsed his account of argument, not his philosophy of mathematics. 
They have typically treated (at least some) mathematical arguments as substantial rather than analytic. For example, ``It is the substantial argumentation that is seen here as more adequate for the reconstruction'' of mathematics classroom situations \cite[236]{Krummheuer95}.
In other words, the part of Toulmin's work that we should apply to many mathematical arguments is what he has to say about 
non-mathematical arguments.

Much subsequent work applying the Toulmin layout to mathematical reasoning has concerned ways in which it may be extended to cover 
cases that he does not directly address---including many of the less favoured sides of the distinctions drawn in Section~\ref{section:Argument} above.
Toulmin's own later work discusses how more than one layout may be chained together to represent a multi-step argument \cite[79]{Toulmin79}.
Other authors have shown how more complicated structures, including linked and convergent arguments, may be represented by combinations of layouts \cite[214]{Aberdein06c}.
Christine Knipping and David Reid have paid particular attention to larger scale structures of parallel argument, distinguishing source-like argumentation structure, where ``arguments and ideas arise from a variety of origins, like water welling up from many springs'', from ``reservoir structure'', wherein arguments ``flow towards intermediate target-conclusions that structure the whole argumentation into parts that are distinct and self-contained'', and spiral argumentation structure in which ``the final conclusion is proven in many ways''
\cite[18 ff.]{Knipping19}. 
And Matthew Inglis and colleagues have argued persuasively for the relevance of the full Toulmin model to mathematical arguments, rather than the simplified DWC or DWBC versions 
that have found widest application amongst mathematics educationists \cite{Inglis0+}. In particular, they observe that non-deductive warrants can play an essential role in mathematical argumentation, 
just so long as this is signalled by the use of appropriate qualifiers.
Work has also been done to explore the connections between the Toulmin layout and other models of argument applicable to mathematics \cite{Pease11};
or to link it to broader conceptual analyses of mathematical cognition, such as the 
``ck\textcent-enriched'' Toulmin model
\cite{Pedemonte16}.

\subsection{Argumentation schemes}
Argumentation schemes are stereotypical patterns of reasoning. Although their origins lie in 
the topoi of classical rhetoric, they have lately found extensive application 
in the analysis and evaluation of argumentation. 
This revival is substantially due to the work of the argumentation theorist Douglas Walton. 
Most attention has been paid to defeasible schemes typical of informal reasoning, although deductive inference rules can also be considered special cases of argumentation schemes.
The defeasible nature of the reasoning is not made explicit amongst the premisses, but captured by an additional device, {critical questions}, which point to possible exceptions. 
Many of the defeasible schemes may ultimately be understood as more or less specialized instances of the very general scheme of \textit{Defeasible Modus Ponens} \cite[366]{Walton08}.
In Scheme~\ref{DMP} we have presented it in a way designed to bring out its similarities to the Toulmin layout:

\begin{scheme}\label{DMP}
Defeasible Modus Ponens
\end{scheme}
\begin{description}
\item[\em Data]$P$.
\item[\em Warrant] As a rule, if $P$, then $Q$.\\
Therefore, \dots 
\item[\em Qualifier] presumably, \dots
\item[\em Conclusion] \dots\ $Q$. 
\end{description}
\begin{center}
Critical Questions
\end{center}
\begin{enumerate}
\item \emph{Backing:} What reason is there to accept that, as a rule, if $P$, then $Q$?
\item \emph{Rebuttal:} Is the present case an exception to the rule that if $P$, then $Q$?
\end{enumerate}

The strength of the argumentation scheme approach lies in its heterogeneity: an influential (but not exhaustive) survey identifies over one hundred different schemes \cite[308 ff.]{Walton08}. 
Hence schemes are typically presented with much greater specificity than Scheme~\ref{DMP}.
For example, Scheme~\ref{Analogy} is a scheme for Argument from Analogy:

\begin{scheme}\label{Analogy}
Argument from Analogy
\end{scheme}
\begin{description}
\item[\em Similarity Premise] Generally, case $C_{1}$ is similar to case $C_{2}$.
\item[\em Base Premise] $A$ is true (false) in case $C_{1}$.
\item[\em Conclusion] $A$ is  true (false) in case $C_{2}$.
\end{description}
\begin{center}
Critical Questions
\end{center}
\begin{enumerate}
\item Are there differences between $C_{1}$ and $C_{2}$ that would tend to undermine the force of the similarity cited?
\item Is $A$ true (false) in $C_{1}$?
\item Is there some other case $C_{3}$ that is also similar to $C_{1}$, but in which $A$ is false (true)? \cite[315]{Walton08}
\end{enumerate}

Analogies in mathematics can be formal and thereby capable of rigorous proof \cite[for a specific example and further discussion of this scheme, see][373]{Aberdein13a}. 
They can also be informal heuristics, for example the 
``strong analogy between the pluralist nature of set theory and what has emerged as an established plurality in the foundations of geometry'' (\citealp[296]{Hamkins20}; for further discussion of this analogy, see \citealp{Berry20}).
Analogical reasoning has been a topic of wide interest and 
mathematical analogies in particular have been the subject of focussed discussion \cite{Schlimm08,Bartha12,Priestley13}.

\setlength{\rotFPtop}{\rotFPbot}
\begin{sidewaystable}
\begin{center}\footnotesize
\begin{tabular}{*{3}{p{.22\textwidth}}p{.25\textwidth}}
\toprule
Discovery arguments&Applying rules to cases&Practical reasoning&Source-dependent arguments\\
\midrule
\begin{enumerate}[leftmargin=1.5em]
\item Arguments establishing rules
\begin{itemize}
\item Argument from a random sample to a population
\item Argument from best explanation
\end{itemize}
\item Arguments finding entities
\begin{itemize}
\item \emph{Argument from sign}\textsuperscript{\emph{f,h,i,k}} 
\item\raggedright \emph{Argument from ignorance}\textsuperscript{\emph{c}} 
\end{itemize}
\end{enumerate}\rule{0pt}{12em}\quad
\fbox{\parbox{.19\textwidth}{\scriptsize
\begin{enumerate}[label=\emph{\alph*}]
\item\cite{Aberdein07b}
\item\cite{Aberdein09}
\item\cite{Aberdein13}
\item\cite{Aberdein13a}
\item\cite{Aberdein19a}
\item\cite{Aberdein21a}
\item\cite{Cantu13}
\item\cite{Dove09}
\item\cite{Metaxas15}
\item\cite{Metaxas16}
\item\cite{Pease11}
\end{enumerate}
}}
&
\begin{enumerate}[leftmargin=1.5em]
\item Arguments based on cases
\begin{itemize}
\item \emph{Argument from an established rule}\textsuperscript{\emph{f}} 
\item\raggedright \emph{Argument from verbal classification}\textsuperscript{\emph{b,i,j}} 
\item Argument from cause to effect
\end{itemize}
\item\raggedright Defeasible rule-based arguments
\begin{itemize}
\item\raggedright \emph{Argument from example}\textsuperscript{\emph{c,e,i,j}}
\item \emph{Argument from analogy}\textsuperscript{\emph{d,j}}
\item \emph{Argument from precedent}\textsuperscript{\emph{k}}
\end{itemize}
\item Chained arguments connecting rules and cases
\begin{itemize}
\item\raggedright \emph{Argument from gradualism}\textsuperscript{\emph{c}}
\item Precedent slippery slope argument
\item Sorites slippery slope argument
\end{itemize}
\end{enumerate}&
\begin{enumerate}[leftmargin=1.5em]
\item\raggedright Instrumental argument from practical reasoning
\begin{itemize}
\item Argument from action to motive
\end{itemize}
\item Argument from values
\begin{itemize}
\item Argument from fairness
\end{itemize}
\item \emph{Value-based argument from practical reasoning}\textsuperscript{\emph{g}}
\begin{enumerate}
\item \emph{Argument from positive or negative consequences}\textsuperscript{\emph{c,i,j}}
\begin{itemize}
\item Argument from waste
\item Argument from threat
\item Argument from sunk costs
\end{itemize}
\end{enumerate}
\end{enumerate}
&
\begin{enumerate}[leftmargin=1.5em]
\item Arguments from position to know 
\begin{enumerate}
\item\raggedright \emph{Argument from expert opinion}\textsuperscript{\emph{a,c,i}}
\item \emph{Argument from position to know}\textsuperscript{\emph{c}}
\begin{itemize}
\item Argument from witness testimony
\end{itemize}
\end{enumerate}
\item Ad hominem arguments
\begin{enumerate}
\item Direct ad hominem 
\item\raggedright Circumstantial ad hominem 
\begin{itemize}
\item Argument from inconsistent commitment
\item Arguments attacking personal credibility
\begin{enumerate}
\item Arguments from allegation of bias
\item Poisoning the well by alleging group bias
\end{enumerate}
\end{itemize}
\end{enumerate}
\item Arguments from popular acceptance
\begin{itemize}
\item\raggedright \emph{Argument from popular opinion}\textsuperscript{\emph{b}}
\item \emph{Argument from popular practice}\textsuperscript{\emph{b}}
\end{itemize}
\end{enumerate}
\\
\bottomrule
\end{tabular}
\end{center}
\caption{Walton \& Macagno's 
partial classification of schemes (adapted from \citealp[22
]{Walton16a}), with 
prior applications to mathematical argumentation indicated.}
\label{schemes}
\end{sidewaystable}%

In more recent work, Walton proposed a partial taxonomy of schemes \cite[22]{Walton16a}, although he acknowledged that some schemes remained outside this classification.
Many of these schemes have been applied to mathematical arguments: Table~\ref{schemes} is based on \cite[22, Table~1]{Walton16a}, but adds references to prior work in which such applications have been developed. 
Just as Walton's classification of schemes is incomplete, so is their application to mathematics. 
Some of these schemes may be of limited usefulness in the analysis of specifically mathematical argumentation, but others have direct application.

By varying which schemes are treated as admissible, it is possible to capture different conceptions of mathematical rigour.
To this end, one of us has proposed a threefold distinction among the ways schemes may be employed in mathematical reasoning \cite[366 f.]{Aberdein13a}:
\begin{itemize}
\item {\bf A-schemes} {correspond directly to derivation rules}. 
(Equivalently, we could think in terms of a single A-scheme, the `pointing scheme' which picks out a  derivation whose premisses and conclusion are formal counterparts of its data and claim.) 
\item {\bf B-schemes} are {exclusively mathematical arguments}: high-level algorithms or macros. Their instantiations {correspond to substructures of derivations} rather than individual derivations (and they may appeal to additional formally verified propositions). 
\item {\bf C-schemes} are even looser in their relationship to derivations, since {the link between their data and claim need not be deductive}. 
Specific instantiations may still correspond to derivations, but there will be no guarantee that this is so and no procedure that will always yield the required structure even when it exists.
Thus, where the qualifier of A- and B-schemes will always indicate deductive certainty, the qualifiers of C-schemes may exhibit more diversity. Indeed, different instantiations of the same scheme may have different qualifiers.
\end{itemize}

B-schemes are essentially what Saunders Mac~Lane calls ``processes of proof'' or 
general rules, algorithmic procedures that are ultimately reducible to elementary logical inferences although not necessarily so analysed by the mathematicians who routinely employ them: 
``whenever a group of elementary processes of proof occurs repeatedly in the course of many proofs, it is desirable to formulate this group of steps once for all as a new process" \cite[123]{MacLane35}.
Much more recently, Yacin Hamami has used B-schemes, which he terms ``hl-rules'' or higher-level rules of inference, to defend 
the ``standard view'' of mathematical rigour, that rigorous proofs are those for which there is a routine translation into a formal derivation \cite{Hamami19}.

Hamami's account of rigour has three components: a descriptive thesis; a normative thesis; and a philosophical thesis asserting the conformity of the other two theses. Relative to some mathematical practice $\mathcal{M}$, these theses may be stated as follows.
The descriptive thesis states that a mathematical proof $P$ is rigorous$_D$ if and only if for every mathematical inference $I$ in $P$, there exist $D \in \mathcal{D}^*$ and $V_1,\dots, V_n \in\mathcal{V}^*$ such that (1) $D(I) = \langle I_1,\dots,I_n\rangle$ and (2) $Vi(Ii) = \mathsf{valid}$ for all $i \in \mathbb{\Lbrack}1,n\mathbb{\Rbrack}$, where
the set $\mathcal{D}^*$ consists of decomposition (or proof search) processes, whereby a mathematical inference may be rewritten as a sequence of immediate mathematical inferences, and the set $\mathcal{V}^*$ consists of hl-rules, whereby immediate mathematical inferences are judged valid if they correspond to instances of the hl-rules.
The normative thesis states that a mathematical proof $P$ is rigorous$_N$ if and only if $P$ can be routinely translated into a formal proof.
Hamami defines ``routine translation'' 
as the composition of three successive translations between proofs understood at four levels of granularity: 
vernacular level proofs, comprised of inferences presented at the level of formality normal to $\mathcal{M}$;
higher-level proofs, comprised of inferences instantiating hl-rules in $\mathcal{M}$;
intermediate-level proofs, comprised of inferences instantiating primitive rules of inference in $\mathcal{M}$; and
lower-level proofs, comprised of inferences instantiating rules of inference in a purely formal system.
The conformity thesis states that if $P$ is rigorous$_D$ then $P$ is rigorous$_N$.
Hamami's account of rigour corresponds to one of four alternatives that one of us has discussed elsewhere \cite[369]{Aberdein13a}. It may be contrasted with the more conservative policy of only admitting A-schemes and thereby treating only formalized mathematics as truly rigorous and more liberal options in which C-schemes are also admitted, thereby tolerating a greater diversity of innovation in mathematical proof.

\section{Mathematical Reasoning as Argumentative}
\label{section:Reasoning}
In Section \ref{section:Proof}, we saw that to understand proof through the lens of argumentation theory is to see it in the context of a greater diversity of types of mathematical argument. While proofs are a vital component of modern mathematical practice, 
it is not the only aspect of mathematical reasoning to which argumentation theory may be applied. 
We now turn to these other aspects. Our discussion begins with general claims about mathematical reasoning and rhetoric. We then turn to a number of issues which surround proof including problem choice, reasoning about refutations, and presentation of mathematical information. 

\subsection{Mathematics \& Rhetoric}
A first connection between mathematical reasoning and argumentation theory involves rhetoric. Rhetoric, the study of the art of persuasive argument, has long been set in opposition to mathematics. It was thought that mathematics, as an objective, rational, and atemporal field, has little to do with 
the study of persuasion. But a number of authors have challenged this idea. An early paper by Philip J.~Davis and Reuben Hersh 
identified two areas where mathematics and rhetoric intersect \cite{Hersh86}. The first 
involves importing or applying mathematics to theories in the social sciences, such as psychodynamics and economics. Such appeals to mathematiziation are rhetorical and argumentative moves, but they are not argumentation within mathematical reasoning.
However, as 
Davis and Hersh point out, there is rhetoric in mathematics too, since ``all proofs are incomplete, from the viewpoint of formal logic'' \cite[66]{Hersh86}. Each proof requires rhetorical elements to convince the intended audience that the result is true. The mathematician relies on the audience's background knowledge or intuition to patch up gaps in the proof and understand the intentions of the prover.  

Like Davis and Hersh, Edward Schiappa 
discusses multiple ways in which mathematical reasoning can be rhetorical. The first intersection, again, is the rhetorical use of mathematics. Mathematical methods can be used to 
persuade in a variety of arguments. 
Schiappa cites examples ranging from mundane advertisement---``four out of five dentists agree''---to the discovery of Neptune to argue that mathematical reasoning plays a role in lending credibility to arguments outside its purview \cite{Schiappa21}. 
The second 
intersection is the role of rhetoric within mathematics: 
the argumentative and stylistic modes of persuasion in mathematical arguments. Each aspect of mathematical practice historically required a social and persuasive component. Acceptance of axioms and stipulated definitions depends on the audience one aims to persuade. Even the available concepts which mathematicians reason about can be the result of rhetoric. For example, G.~Mitchell Reyes argues that, since there was no empirical or geometric verification for infinitesimals, their substance was found in the rhetorical argumentation which surrounded them \cite{Reyes04}. Schiappa also argues that the language of mathematics is rhetorical since it is human-made. Symbols and concepts like the infinitesimal or the number zero were additions where ``social acceptance was not assured, meaning was contested, and alternatives competed'' \cite[49]{Schiappa21}. Relatedly, in this collection, Reyes examines the relationship between rhetoric and mathematics in Lakatos's \textit{Proofs and Refutations} \cite{Reyes21}.  

\subsection{Problem Choice}

In addition to the rhetorical components pervasive in non-proof, there is room for argumentation in the process around proofs. Perhaps the most fundamental part of solving a problem, and of proving a theorem, is selecting a problem. Problem choice in mathematics has frequently been associated with the beauty and intrinsic worth of problems. Under this view, there is 
a special sensibility mathematicians employ when choosing a subject. According to Jacques Hadamard, mathematicians ``feel that such a direction of investigation is worth following; [they] feel that the question in itself deserves interest \dots\ everybody is free to call or not to call that a feeling of beauty'' \cite[127]{Hadamard45}. This approach to problem choice has been supplemented in recent years.  
Morten Misfeldt and Mikkel Willum Johansen interviewed research mathematicians about the factors which influence problem choice \cite{Misfeldt15}. In line with Hadamard's claims, certain external factors, like funding, were not very influential on problem choice. 
Misfeldt and Johansen found that problem choice was largely motivated by personal interest, perceived ability to solve the problem, and the values of the community. It was not enough for mathematicians to be interested in the problem, they had to be assured that the mathematical community would be interested. Mathematicians expressed ``the need to have an audience---the right audience---for their work'' \cite[368]{Misfeldt15}. This connection between an audience and an arguer is an area which argumentation is primed to explore. 

Elsewhere, one of us has looked at problem choice through an argumentative lens by 
applying 
Cha\"im Perelman and Lucie Olbrechts-Tyteca's notion of the contact of minds 
to problem choice in mathematics \cite{Ashton18}. The contact of minds is a set of conditions which must be met before argumentation can occur. Contact of minds requires four things: the arguer must attach importance to the audience, the speaker must not be beyond question, the audience must be willing to consider being convinced, and they must share a common language. Contact of minds is a prerequisite of any argumentation and mathematical arguments are no different. Ashton argues that it is an audience-based factor that influences problem choice alongside traditional considerations of beauty, intrinsic worth, and practical benefits. 
But the contact of minds was not originally meant to be applied to mathematics: 
Perelman and Olbrechts-Tyteca specifically oppose 
their study of argumentation 
to the mathematical sciences. Argumentation, they claimed, was distinct from mathematics since argumentation was both social and the conclusions were probable \cite{Perelman69}. But barring mathematical practice from the domain of argumentation is inapproporiately limiting \cite{Dufour12, Ashton18, Ashton20}. Misfeldt and Johansen's interviews indicate that mathematicians consider the interest of other mathematicians while choosing problems. Given that problem choice is broader than merely the structure of proof, a choice about what to research is a choice of what to argue about. In this way, 
the process of choosing a problem to research 
is distinctly social and related to the mathematical community 
(\citealp{Ashton18}; discussed further in Section~\ref{section:Communities}). 

\subsection{Reasoning about Refutations}\label{subsection:Refutations}
The influential informal logician Ralph Johnson \citeyearpar{Johnson00} 
follows the rhetoricians Perelman and Olbrechts-Tyteca \citeyearpar{Perelman69} 
in denying that proofs can be arguments because there are features of proof, 
such as necessity, that are incompatible with the social dimension of arguments. For Perelman and Olbrechts-Tyteca, mathematics deals in demonstrations which are mechanically checkable and result in certainty, regardless of the audience. Mathematical reasoning, for them, does not involve uncertainty or controversy. Likewise, Johnson viewed proofs as lacking a dialectical tier. For Johnson, arguments have two components. The first component of argument is the illative core which is a discursive structure where reasons support the conclusion. But this logical aspect alone is not a full argument. A dialectical tier must supplement the illative core. Dialectical tiers are where arguers discharge their dialectical obligations by responding to objections, criticisms, or implications of their view. According to Johnson, proofs don't have dialectical tiers since they are conclusive and anyone trained in the field must recognize that they are conclusive. But refutations are a natural part of mathematical reasoning. Reasoning about purported refutations may best be understood in terms of argumentation.

In an idealized view of proof, one in which proof cannot be an argument, the reader of a proof is an expert who ``needs nobody to grasp a proof, otherwise she is not an expert'' \cite[69]{Dufour12}. Such a view, according to Dufour, ignores the important role of checking a proof for correctness. There may still be legitimate room for refutations at this stage. As \citet{Fallis03} notes, proofs have many intentional gaps. The gaps are purportedly something a reader could fill in, with enough time and background knowledge. Verifying that a proof is correct requires checking these gaps. If the gaps cannot be traversed by either the reader or author, the proof itself may be refuted. Sometimes large gaps may result in the apparent refutation of otherwise good proofs. According to Dufour, Galois experienced such problems of communication \cite[71]{Dufour12}. Galois's exaggerated brevity led his audience to believe that the mathematics itself was insufficiently developed. Proofs may be refuted because of unintelligible or untraversable gaps. These gaps, and their relationships to audience understanding, are evidence that reasoning about the correctness of proofs, and their refutations, involves argumentation. 

Contra Johnson, Ian Dove has argued that proofs do have dialectical tiers \cite{Dove07}. In particular, Dove argues that the method of monster barring seen in Lakatos's \textit{Proofs and Refutations} is part of a dialectical tier \cite{Lakatos76}. 
Lakatos reconstructs a series of purported proofs of Euler's formula. Cauchy's proof of Euler's formula for polyhedra is considered and then counterexamples are raised which are not convex and not simply connected. Dove argues that there is a dialectical tier for Cauchy's proof of Euler's formula since (a) objections are raised to the proof and (b) the objections receive responses in the proof. The method of response is to bar exceptions to the formula so that Euler's formula holds for simply connected, convex polyhedra. In other words, monster barring to improve a conjecture is one way a proof can have a dialectical tier. 

Both Dove and Dufour found that reasoning about refutations is an argumentative practice. Purported proofs are not always uncontroversial and incorporating refutations involves filling out a dialectical tier or considering a relevant audience. 

\subsection{Presenting Reasoning}
After a problem has been chosen and its solution has been verified, mathematical reasoning must be disseminated to a larger public. There is much work to be done to discover how mathematicians incorporate audience consideration while presenting their solutions. But Line Edslev Andersen's interviews with working mathematicians indicate that audience consideration does feature into how papers are revised. Andersen's interviews provide insight into how mathematicians write for mathematicians \cite{Andersen20a}, how peer reviewers receive and evaluate papers \cite{Andersen17}, and how mathematicians choose which gaps to include in their papers \cite{Andersen20}. Additionally, mathematical results must be translated to students in a pedagogical setting. Results also often need to be communicated to scientists or mathematicians in other fields. In each of these cases, it is common to reword or recast a proof to aid in communication. All ``strategic rewordings belong to the field of mathematical argumentation'' \cite[73]{Dufour12}. Some of these rewordings may even be usefully cast within Johnson's concept of a dialectical tier. This is a rich area for further research, but, as Dufour points out, one must first admit mathematical proofs into the realm of argumentation. Section \ref{section:Communities} 
examines how communities play a role in mathematical practice broadly, not just in presenting reasoning. 

\section{Mathematical Communities as Argumentative}
\label{section:Communities}
Throughout this chapter we have looked at the different aspects of mathematical practice that involve argumentation. Mathematicians present and produce argumentations. Their choice of problems and reasoning while solving problems also involved argumentation. But arguments are presented to, or developed with, an audience in mind. We now turn to mathematical audiences and communities. We argue that mathematical audiences influence which investigations are undertaken and how they are undertaken. 

Given that the interest of this section is mathematical audiences, the receivers of mathematical arguments, one might begin by asking what kind of interaction mathematicians have with their audiences. Mathematics is often portrayed as an isolated activity. The mathematician locks herself in an attic and spends days proving complicated conjectures alone. Indeed, stories of famous mathematicians seem to support such an idea. Andrew Wiles did work on his proof of Fermat's Last Theorem in an attic. And mathematical advances can be so particular, as in the case of 
Mochizuki's purported proof of the $abc$ conjecture, that only a handful of people can verify or understand them. But this isolated view is not an accurate portrayal of mathematical practice.  

Without an appropriate community, or engagement with that community, mathematical investigations can falter. Besides the obvious importance of mentors or co-authors, there needs to be a certain level of engagement with a larger community of mathematicians. Take, for example, William Thurston's discussion in ``On Proof and Progress in Mathematics'' \cite{Thurston94}. According to Thurston, mathematicians tend to follow fads. Fad-following is in line with his claim that there is a vital social component to mathematical progress. He drew this conclusion from his own experiences. When Thurston first entered the field of foliations, he rapidly solved a number of open problems. But the field seemed to empty out within a few years of his entrance and Thurston lost interest soon after. In response to this experience, he tried to develop and present infrastructure instead of theorem-proving. He writes that:
\begin{quote}
	I have put a lot of effort into non-credit-producing activities that I value just as I value proving theorems: mathematical politics, revision of my notes into a book with a high standard of communication, exploration of computing in mathematics, mathematical education, development of new forms for communication of mathematics \dots\ I think that what I have done has not maximized my ``credits'' \dots\ I do think that my actions have done well in stimulating mathematics. \cite[177]{Thurston94}
\end{quote}
We can see that, for Thurston, the health of the mathematical community is invaluable to retaining mathematicians who are interested in that area. This is part of their tendency to `follow fads' to newly exciting communities. 

Audience and community involvement seem important, but 
what exactly does argumentation theory have to offer in this area? As one of us has argued, Thurston's story exemplifies a broader issue namely that problem choice in mathematics rests on the assurance that there exists a `contact of minds' between the audience and the mathematician \cite{Ashton18}. As we can see, the interest and activity of a mathematical community is vital to successful mathematical practice. 

In addition to the role that communities play in problem choice, they play an important role in the ongoing dialectic involved in problem solving. The method of proofs and refutations described by Lakatos relies heavily on community involvement in problem-solving. A first `proof' is suggested to show that, for polyhedra, $V - E + F = 2$, where $V$ is the number of vertices, $E$ the number of edges, and $F$ the number of faces of the polyhedron. The proof faces a number of global and local counterexamples, that is counterexamples to the conjecture and to specific proof steps, respectively. Accommodating these counterexamples requires a reconceptualization of the definition of polyhedron. In Section~\ref{section:Reasoning}, we saw that reasoning about those refutations involved filling out a dialectical tier. It is also important that community involvement plays an essential role in the reasoning. The history of the proof is recounted using the story of a class attempting to solve the problem. But it is clear that the refinements of definition, the counterexamples, and the resulting proofs are all products of dialectical engagement between the initial prover and an audience. 

This community involvement plays out on a historical scale, as in \cite{Lakatos76}, but also in other collaborative problem-solving. For example, communities play a key role in {crowdsourced} mathematics such as {Mathoverflow} and {(Mini-)Polymath}. These activities have been studied by Ursula Martin and Alison Pease who connect the activities undertaken to the method of proofs and refutations \cite{Martin13a, Martin13b}. Building on these data sets, in later work with Joseph Corneli and other colleagues, they model mathematical arguments by analyzing how the discourse unfolds \cite{Corneli18}. They use what they call Inference Anchoring Theory + Content to further understand how these dialogues introduce and track salient features of mathematical progress. This is done by examining the dialogue within a community of mathematicians. 

So far we have considered the role that an external audience plays in problem choice. We have also looked at how we could identify communities in terms of the argument schemes that they allow or even disavow. But the final role of mathematical audiences is best described as an internal one. The core idea is that mathematical audiences play an important role in the development of proofs and other mathematical arguments. Valentin Bazhanov, for example, has argued that proofs are an appeal to a scientific community \cite{Bazhanov12}. A new result is offered to the community as a purported proof. To become a proof, the community must be persuaded that the argument is reliable and reproducible. This decision is borne out in dialogue. 

Catarina Dutilh Novaes argues for a dialogical conception of proof \cite{DutilhNovaes21}. 
According to the dialogical conception, the concept of proof is an embodiment of a semi-adversarial dialogue between two people: Prover and Sceptic. Sceptic grants certain premises to Prover. Prover then puts forward statements claimed to follow from the premises. At each move, or inference, Sceptic has three potential moves. He brings up objections, counterexamples, and requests for clarification. If all of the steps are indefeasible, that is without counterexample, the proof is a winning strategy for Prover. Of course proofs are not usually dialogues between two interlocutors following Prover-Sceptic rules. Dutilh Novaes argues that the Sceptic is actually internalized into the method itself. In this sense, Sceptic, who is a particular kind of audience, plays a vital role in what constitutes proof. 

In addition to the internalized sceptic, one of us has argued that the audience under consideration in a proof is actually a universal audience \cite{Ashton20}. 
This is a concept drawn from Perelman and Olbrechts-Tyteca: an argument to the universal audience is one that is meant to convince all people. The universal audience itself is an imagined audience---an arguer can never stand before all people and ask whether or not they assent to his argument. Rather, the universal audience is an abstraction from experiences with real audiences that the arguer has encountered. The account applies also to mathematics: mathematicians encounter real audiences in their education and practice. In addition, they learn what certain groups of real audiences react to---knot theorists accept inferential moves involving diagrams that certain algebraists might not. By considering all these different, real audiences, the mathematician constructs an internalized audience which reflects the standards of reasonableness found within each audience. In this way, the mathematician constructs his own universal audience. 

According to both Dutilh Novaes and Ashton, 
there is an internalized audience which, through experience with real audiences, is vital to proof development. In one case, a proof is an argument to an internalized sceptic. In the other, it's an argument to an internalized standard of reasonableness. Whether the mathematician aims to convince the adversarial sceptic or `all reasonable people,' there is clearly a role for the audience in standards of proof. In other words, proofs should be examined in terms of the audiences which could have influenced them.

Communities, conceived of as audiences to mathematical arguments, are a vital component of the practice on a number of levels. The assurance of an interested audience and the contact of minds can influence research programmes and problem choice. The ongoing dialogue between provers and the mathematical community helps to verify results, introduce new methods, and clear hidden assumptions. In addition to these explicit, external roles, the internalized audience plays an important role in proof development and the associated normativity. Given all this, it's clear that no account of argumentation in the philosophy of mathematics could be complete without a thorough discussion of the role of these audiences. 


\end{document}